\input amstex
\documentstyle {amsppt}
\magnification=1200
\vsize=9.5truein
\hsize=6.5truein
\nopagenumbers
\nologo

\def\norm#1{\left\Vert #1 \right\Vert}
\def\abs#1{\left\vert #1 \right\vert}
\def\Mod{\text{\rm{Mod}}}
\def\PMF{\Cal{PMF}}
\def\G{\Gamma}
\def\On{\text{\rm{Out}}(F_n)}
\topmatter

\title
Mapping~class~groups
and~outer~automorphism~groups~of~free~groups
are $C^*$-simple
\endtitle

\rightheadtext{Mapping class groups and $\text{\rm{Out}}(F_n)$}

\author Martin R.~Bridson and Pierre de la Harpe
\endauthor

\thanks
The authors acknowledge support from the
{\it Swiss National Science Foundation}.
\endthanks


\endsubjclass

\abstract
We prove that the reduced $C^*$-algebras of centerless
mapping class groups and  outer automorphism
groups of free groups are simple,
as are the irreducible pure subgroups of mapping class groups and the
analogous subgroups of outer automorphism groups of free groups.
\endabstract

\address
Martin Bridson, Mathematics Department, Huxley Building,
\newline Imperial College
London, SW7 2AZ, UK.
\newline
E-mail: m.bridson\@imperial.ac.uk
\endaddress

\address
Pierre de la Harpe, Section de Math\'ematiques, Universit\'e
de Gen\`eve, C.P. 240,
\newline
CH-1211 Gen\`eve 24, Suisse.
\newline
E-mail: Pierre.delaHarpe\@math.unige.ch
\endaddress
 \endtopmatter

\document

\head
Introduction
\endhead
\medskip

A group
$\Gamma$ is {\it $C^*$-simple} if its
{\it reduced $C^*$-algebra} $C_\lambda^*(\Gamma)$
 is simple as a complex algebra (i.e.
has no proper two-sided ideals).
The main purpose of this note is to
explain why mapping class groups of surfaces of finite type and
outer automorphism groups of free groups are $C^*$-simple.

There is a fascinating analogy between lattices in semi-simple
Lie groups, on the one hand, and mapping class groups and outer
automorphism groups of free groups, on the other. The results
recorded here resonate well with this analogy, as we shall
explain in a moment.
First, though, we remind the reader that the $C^*$-simplicity of
a group $\G$ may be regarded as a
property of the unitary representation theory of the group. By
definition, $C^*_\lambda(\G)$ is the norm closure of the image
of the complex group algebra $\Bbb C [\G]$ under the left-regular
representation $\lambda_\G : \Bbb C [\G]\to \Cal L (\ell^2(\Gamma))$
defined for $\gamma\in\G$
by $\left(\lambda_{\Gamma}(\gamma)\xi\right)(x) =
\xi(\gamma^{-1}x)$
for all $x \in \Gamma$ and $\xi \in \ell^2(\Gamma)$. A group
$\Gamma$ is $C^*$-simple if and only if
any unitary representation $\pi$ of $\Gamma$ which is weakly contained in
$\lambda_{\Gamma}$ is weakly equivalent to $\lambda_{\Gamma}$
(see Theorem 3.4.4 and Proposition 18.1.4 in \cite{D$C^*$--69}).
Examples of $C^*$-simple groups include
non-abelian free groups \cite{Pow--75},
non-trivial free products \cite{PaS--79},
torsion-free, non-elementary hyperbolic groups
(see \cite{Har--85} and \cite{Har--88}),
and Zariski-dense subgroups in centerless, connected, semisimple
real Lie groups with no compact factors \cite{BCH2--94}.
(In particular $\text{\rm{PSL}}(n, \Bbb Z)$ is
$C^*$-simple \cite{BCH1--94}.)

Apart from the usual low-genus exceptions, mapping class groups and outer
automorphism groups of free groups do not lie in any of
the above classes of groups. But
whenever one has an interesting property of groups with such a
list of examples, the analogy between lattices and
these important groups demands attention. Previous experience
encourages us with many examples where the analogy goes through,
and even leads us to expect that the shape of
proof used in the classical setting might
thrive in the transplanted environment.
This is the case, for example, with theorems concerning
aspects of
rigidity \cite{BrV--00}, \cite{BrV--01},
\cite{FaM--98}, \cite{Iva--97}
and homological stability \cite{Hat--95},  \cite{Ha--85}.
It is also the case for analogues of the Tits Alternative
\cite{McC--85}, \cite{Iva--84}, \cite{BLM--83},
\cite{BFH--00},  \cite{BFH2},  \cite{BFH3}. See also \cite{Vog--02},
\cite{Iva--92} and \cite{Bes--02}.

As well as providing us with encouragement,
these works
also warn us that the tools required to adapt proofs
from the classical setting to that of mapping class
groups and outer automorphism groups of free groups are often
 non-trivial and may require
significant innovation. Fortunately, though, in the
present setting, the necessary tools can be readily gleaned from
work of previous authors, as we shall now explain.

Although mapping class groups and outer automorphism groups of
free groups are more often thought of in analogy with
higher rank lattices (with $\text{\rm{SL}}(n,\Bbb Z)$
springing most readily to mind), it is widely recognized
that there are a number of important respects in which they
exhibit rank-one phenomena  (see \cite{FLM--01}). The argument
that we shall present here falls into the latter category.
Indeed this is what makes our project straighforward:
the argument used in \cite{BCH-94} (following \cite{HaJ--81})
 to establish the $C^*$-simplicity
of rank one lattices (also hyperbolic groups and non-trivial
free products) rests on an elementary lemma, the input of which
is dynamical information reminiscent of the classical ping-pong
lemma, and  the output of which allows one to follow
Powers' proof of the $C^*$-simplicity
of non-abelian free groups. In the case of rank one lattices (more
generally Zariski-dense subgroups) the space on which one
studies the dynamics is the boundary of the symmetric space. At the heart
of the analogy between lattices and mapping class groups lies
the fact that the Teichm\"uller space plays the r\^ole of the
symmetric space;  for outer automorphism groups of free
groups, the r\^ole of
the symmetric space is played by  Culler and Vogtmann's Outer Space
 \cite{CuV--86}.

The mention of ping-pong in this setting
brings to mind the construction of free groups in the
proof of the Tits alternative \cite{Tits--72}, \cite{Har--83}.
The Tits alternative has been established for both
mapping class groups and outer automorphism groups of free groups.
In the case of outer automorphism groups of free groups, this is
recent work that represents the culmination of
a long-term project by Bestvina, Feighn and Handel
developing \lq\lq train-track technology\rq\rq \ \cite{BeH--92}
to produce  refined
topological representatives of free group automorphisms. This
project was motivated by Thurston's work on train-track
representatives of surface automorphisms, which is closely related
to his work on the boundary of Teichm\"uller space. It is the
action of the mapping class group on this boundary that provides
us with the ping-pong behaviour required to establish $C^*$-simplicity
for the mapping class group. In the case of the (outer) automorphism
groups of free groups, we appeal to the work of Bestvina-Feighn-Handel
\cite{BFH--97}
and a refinement of Levitt and Lustig \cite{LeL--03} concerning the dynamics
of automorphisms on the boundary of Outer Space.

\medskip

   Our account of the action of $\text{\rm{Out}}(F_n)$ on the
boundary of Outer Space owes a great deal to the insights of
Karen Vogtmann. We are most grateful to her
for sharing these and other insights during conversations
in Geneva in the summer of 2002.

\bigskip
\head
1. Powers' Criterion for $C^*$-Simplicity
\endhead

We define a  group $\Gamma$ to be  a {\it Powers group} if
\vskip .3truecm
\hskip.5truecm{for any finite subset $F$ in $\Gamma \smallsetminus \{e\}$
and for
   any integer $N \ge 1$,
\par\hskip.5truecm
   there exists a partition  $\Gamma = C \sqcup D$ and elements
   $\gamma_1,\hdots,\gamma_N$ in $\Gamma$    such that
\par\hskip.5truecm
$fC \cap C = \emptyset$ for all $ f \in F$,
  and $\gamma_jD \cap \gamma_kD = \emptyset$
   for all $j \neq k$ in $\{1,\hdots,N\}$.}
\medskip

\smallskip\noindent
The terminology is in honour of  \cite{Pow--75}. The first thing to be said
about this definition is that it implies $C^*$-simplicity;
the argument is essentially that of the original paper \cite{Pow--75}
and is repeated here as an Appendix for the reader's convenience.
The other important thing to note is
that  there is a simple  dynamical criterion that enables one to show
that many interesting groups are Powers groups. In order to describe
this criterion we need the following vocabulary.

A homeomorphism $\gamma$ of a Hausdorff space $\Omega$
is said to be {\it hyperbolic}
if it has two fixed points $s_{\gamma},r_{\gamma} \in \Omega$ and exhibits
{\it north-south dynamics}:
for any pair of neighbourhoods $S$ of $s_{\gamma}$ and $R$ of $r_{\gamma}$,
there exists $n_0 \in \Bbb N$ such that
$\gamma^n(\Omega \smallsetminus S) \subset R$ and
$\gamma^{-n}(\Omega \smallsetminus R) \subset S$
for all $n \ge n_0$.
The points $s_{\gamma}$ and $r_{\gamma}$ are  called the
{\it source} and the {\it range} of $\gamma$, respectively.
 \par

   Two hyperbolic homeomorphisms of $\Omega$ are {\it transverse}
if they have no common fixed point. \par

\proclaim{1.1\ Proposition} Let $\Gamma$ be a group acting by homeomorphisms
on a
Hausdorff space~$\Omega$. Assume that the following two conditions hold.
\roster
\item"(i)"
$\Gamma$ contains two transverse
hyperbolic homeomorphisms of $\Omega$.
\item"(ii)"
For any finite subset $F$ of $\,\Gamma \smallsetminus \{e\}$,
there exists a point $t \in \Omega$ fixed by some hyperbolic homeomorphism
of
$\,\Gamma$ such that $ft \neq t$ for all $f \in F$.
\endroster
  Then $\Gamma$ is a Powers group, and in particular $C^*_{\lambda}(\Gamma)$
is a simple $C^*$-algebra.
\endproclaim

\demo{Proof} Consider $F \subset \Gamma \smallsetminus \{e\}$ and $N \ge 1$
as in the definition of a Powers group.
By hypothesis, there exist hyperbolic homeomorphisms
$\gamma, \gamma', \gamma'' \in \Gamma$
and a neighbourhood $C_{\Omega}$ of the range $r$ of $\gamma$
such that $\gamma',\gamma''$ are transverse
and such that $fC_{\Omega} \cap C_{\Omega} = \emptyset$ for all $f \in F$.
Let $\gamma_1,\hdots,\gamma_N$ be pairwise transverse conjugates of
$\gamma'$
by appropriate powers of $\gamma''$.
Upon conjugating $\gamma_1,\hdots,\gamma_N$ by a large power of $\gamma$,
we may assume that,  for each $j \in \{1,\hdots,N\}$,
both the source $s_j$ and the range $r_j$ of $\gamma_j$ are in $C_{\Omega}$.
We choose neighbourhoods $S_j$ of $s_j$ and $R_j$ of $r_j$ in such a way
that
$S_1, R_1, \hdots, S_N, R_N$ are pairwise disjoint and inside $C_{\Omega}$.
Upon replacing now each of $\gamma_1,\hdots,\gamma_N$ by a large enough
power
of itself, we may furthermore assume that
$\gamma_j(\Omega \smallsetminus C_{\Omega}) \subset R_j$,
and in particular that the $\gamma_j(\Omega \smallsetminus C_{\Omega})$ are
pairwise
disjoint  subsets  of $\Omega$. \par

   Choose $\omega \in \Omega$;
let $C = \left\{ \gamma \in \Gamma \mid \gamma \omega \in C_{\Omega}
\right\}$
and
$D = \left\{ \gamma \in \Gamma \mid \gamma \omega \notin C_{\Omega}
\right\}$.
Then $fC \cap C = \emptyset$,  since
$f C_{\Omega} \cap C_{\Omega} = \emptyset$,  for all $f \in F$,
and $\gamma_j D \cap \gamma_k D = \emptyset$,
since $R_j \cap R_k = \emptyset$, for all $j \neq k$.
$\square$
\enddemo

The following perturbation of Proposition 1.1 is well-adapted to
the examples in which we are interested. We write $\text{\rm{Stab}}(x)$
to denote the stabilizer in $\G$ of a point $x\in\Omega$.

\proclaim{1.2\ Corollary} Let $\G$ be a group acting by homeomorphisms
on a Hausdorff space~$\Omega$.
Assume that the following two conditions hold.
\roster
\item"(i)"
$\Gamma$ contains a hyperbolic
homeomorphism $\gamma_0$ with source $s_0$ and range $r_0$.
\item"(ii)"
There exists a non-trivial element $\gamma_1\in\G$
such that for each integer $i\neq 0$,
the set $$ \gamma_1^{-i} \left( \text{\rm{Stab}}(r_0)  \cup
\text{\rm{Stab}}(s_0)  
\right) \gamma_1^i
\  \cap  \
\left( \text{\rm{Stab}}(r_0)  \cup \text{\rm{Stab}}(s_0)  \right) $$ is just
$\{e\}$.
\endroster
Then $\Gamma$ is a Powers group, and in particular $C^*_{\lambda}(\Gamma)$
is a simple $C^*$-algebra.
\endproclaim

\demo{Proof} The hyperbolic homeomorphisms
$\gamma_0$ and $\gamma_1^{-1}\gamma_0
\gamma_1$ are transverse. Indeed for any positive integer $n$, the
hyperbolic homeomorphisms $\gamma_i:=\gamma_1^{-i}\gamma_0\gamma_1^i,
\ i=1,\dots,n$ are pairwise transverse and no element of $\Gamma$ fixes
a source or range of more than one of them, by condition (ii). We  write
$s_i$ and $r_i$ to denote the source and range of $\gamma_i$.

Given a finite
set $F$, of cardinality $m$ say, we choose $n$ so that $2n>m$. Since
each element of $\Gamma$  fixes at most  two of the $2n$ points
$s_1,r_1,\dots,s_n,r_n$, at least one point on this list is moved
by every element of $F$.
$\square$
\enddemo

\bigskip
\head
2.\ Mapping~Class~Groups~and~Outer~Automorphism~Groups~of~Free
~Groups Satisfy the Powers Criterion
\endhead

Let $\Mod_S$ denote  the group of isotopy classes of orientation-preserving
homeomorphisms
of a compact orientable  surface $S$ (which may have non-empty
boundary). For an introduction to the properties of such groups, see
Ivanov's excellent
survey \cite{Iva--02}.
To avoid the well-known quirks associated with small examples
\footnote{
In each of these exceptional cases $\Mod_S$ has a non-trivial centre,
and hence the reduced $C^*$ algebra is not simple.
However these cases  can be dealt with individually.
For example, if $S$ is a closed surfaces of genus $1$ or $2$,
the quotient of the group  $\Mod_S$
by its centre of order $2$ is a $C^*$-simple group.}  
we assume that $S$ is neither a sphere with $\le 4$ boundary circles,
nor a torus with $\le 2$ boundary circles,
nor a closed surface of genus 2.
In all of the  remaining cases,
the centre of $\Mod_S$ is trivial and the group acts effectively on  the
Thurston boundary of the associated Teichm\"uller space.
There is a natural identification of the Thurston boundary  with  the space
of projective measured foliations
\footnote{Alternatively, laminations.}
$\PMF\!_S$ on $S$ (the basic reference
for this material is  \cite{FLP--79},
while Ivanov's monograph \cite{Iva--92} and survey \cite{Iva-02}
provide an excellent overview).
$\PMF\!_S$ is a topological sphere; in
particular it is Hausdorff. We shall apply the
considerations of Section 2 to the action of $\Mod_S$ on this space.
\par

   A {\it pseudo-Anosov} in $\Mod_S$ is an element that acts as a
hyperbolic homeomorphism of $\PMF\!_S$. It is well-known that $\Mod_S$
contains transverse pairs of pseudo-Anosov classes --- see, for example,
Lemma 2.5 in \cite{McP--89} or Corollary 7.15 in \cite{Iva--92}.

The stabilizers of the fixed points of pseudo-Anosovs are understood
(see Lemma 2.5 in \cite{McP--89} and Lemma 5.10 in
\cite{Iva--92}, for example):

\proclaim{2.1\ Lemma} If $\phi\in\Mod_S$ is pseudo-Anosov, then
the stabilizer of each of its fixed points in $\PMF\!_S$ is virtually
cyclic.
\endproclaim

In outline, one proves this lemma as follows.
The fixed points of a pseudo-Anosov $\phi$ are the projective classes of its
stable and
unstable laminations. If  $\psi\in\Mod_S$ fixes one of these points $[\mu]$,
then it multiplies the measure
on the underlying lamination by a constant factor, $\lambda(\psi)$ say. The
map $\psi \mapsto\lambda(\psi)$
is a homomorphism from the stabilizer of $[\mu]$ in $\Mod_S$ to the
multiplicative group of
positive reals; the image of this homomorphism
is discrete (hence cyclic), the image of $\phi$ is non-trivial, and the
kernel is finite.
It follows that if $\gamma_1$ is a pseudo-Anosov such
that $\phi$ and $\gamma_1$ do not have common powers, then $\phi$ and
$\gamma_1$ are transverse (and hence have powers
that generate a non-abelian free group).
Moreover, for typical (but not all) $\phi$,
the stabilizer of $[\mu]$ is actually cyclic,
generated by $\gamma_0$ say.
Now $\gamma_0$ and $\gamma_1$ satisfy the conditions of Corollary 1.2.  Thus
we have:

\proclaim{2.2\ Theorem}  $\Mod_S$ is a Powers group; in particular its
reduced $C^*$-algebra is simple.
\endproclaim

A subgroup  $\Gamma\subseteq\Mod_S$ is called {\it reducible} if there is a
non-empty closed 1-dimensional submanifold
$C\subset S$ such that for every $f\in\Gamma$ there is a homeomorphism $F$
in the isotopy class $f$ with $F(C)=C$.
The {\it pure} elements of $\Mod_S$ (the definition of which is
somewhat technical)
contain a torsion-free subgroup of finite index in $\Mod_S$.
Every non-trivial
irreducible subgroup  $\Gamma\subseteq\Mod_S$
consisting of pure elements contains a pseudo-Anosov (see Theorem 5.9 of
\cite{Iva--92}). Arguing as above, it
follows easily from Lemma 2.1 that if $\Gamma$ is not cyclic, then it
contains a  pair of transverse pseudo-Anosovs
$\{\gamma_0,\gamma_1\}$ satisfying the conditions of Corollary 1.2.
Moreover, in this context we
do not need to exclude the low-genus exceptions mentioned at the beginning
of this section.

\proclaim{2.3\ Theorem} For every compact surface $\Sigma$,
 every non-cyclic, pure, irreducible subgroup
of $\Mod_\Sigma$ is a Powers group.
\endproclaim

 We now turn our attention to $\On$, the group of outer
 automorphisms of a free group of rank $n$.
 By definition, $\On = \text{\rm{Aut}}(F_n)/\text{\rm{Inn}}(F_n)$,
 where $ \text{\rm{Aut}}(F_n)$ is the  group of automorphisms
of $F_n$ and $\text{\rm{Inn}}(F_n)$ is the group of inner
automorphisms (conjugations).
In the case $n = 2$, the natural map $\text{\rm{Out}}(F_2)
\to\text{\rm{GL}}(2,\Bbb Z)$
is an isomorphism, hence the centre of $\text{\rm{Out}}(F_2)$ has order
two, and the quotient by this centre
is $C^*$-simple (see \cite{BekH--00}).

Henceforth we assume that
$n \ge 3$. In this case it is easy to
check that the centre of $\text{\rm{Out}}(F_n)$ is trivial.
Let $\Cal C$ denote the set of conjugacy classes of $F_n$.
We consider the vector space $\Bbb R ^{\Cal C}$
of real-valued functions on $\Cal C$, equipped with the product
 topology, and the corresponding projective space
$\Bbb P \Bbb R ^{\Cal C}$ with the quotient topology. (Note
that these spaces are Hausdorff.)
The natural action of $\On$ on $\Cal C$ induces an action on
$\Bbb P \Bbb R ^{\Cal C}$.

Culler and Vogtmann's Outer Space $X_n$ may be described as the
space of equivalence classes of free actions of $F_n$
by isometries on $\Bbb R$-trees. Given such an action on a tree $T$,
we associate to each $w\in F_n$ the positive number
$\|w\| = \inf\{d(wx,x)\mid x \in T\}$. The number $\|w\|$ depends only
on the equivalence class of $w$, and the function $w\mapsto \|w\|$
completely determines the equivalence class of the action.
Thus we obtain a natural equivariant injection $j:X_n\hookrightarrow
\Bbb P \Bbb R ^{\Cal C}$. The set $\overline{j(X_n)}\smallsetminus j(X_n)$
is called {\it the boundary of outer space} and is denoted
$\partial_{\infty}(X_n)$. This space is compact.
 For a survey of these and
related matters, see \cite{Vog--02}.

An element $\gamma \in \On$ is an {\it iwip} (irreducible
with irreducible powers) if no proper free factor of $F_n$
is mapped to a conjugate of itself by a non-zero power of any
representative $\tilde \gamma \in \text{\rm{Aut}}(F_n)$ of $\gamma$.
Bestvina, Feighn and Handel show that iwip outer automorphisms of
free groups behave in close analogy with pseudo-Anosov automorphisms
of surfaces.
In particular, in \cite{BFH--97} they define a set $\Cal{IL}$ of
\lq\lq stable laminations\rq\rq \ on which $\On$ acts. Each iwip 
$\phi\in\On$
has two fixed points $\Lambda_\phi^+, \Lambda_\phi^-\in\Cal{IL}$.
Theorem 2.14 of \cite{BFH--97}, whose proof is closely analogous
to that of Lemma 2.1 sketched above,  amounts to the following statement:

\proclaim{2.4\ Lemma} If $\phi\in\On$ is an iwip then the
stabilizers of the fixed points of $\phi$ in $\Cal{IL}$ are virtually
cyclic.
\endproclaim

And as in the case of the mapping class group, one can choose
$\phi$ so that these stabilizers are actually cyclic.
This information about stabilizers of the fixed points of $\phi$
can be transferred to the action of $\On$ on $\partial_\infty(X_n)$
by virtue of
Corollary 3.6 of \cite{BFH--97}:

\proclaim{2.5\ Lemma} There is an $\On$-equivariant injection
$\Cal{IL}\hookrightarrow\partial_\infty(X_n)$.
\endproclaim

As in the proof of  Theorem 2.2 we will be in a position to apply Corollary
1.2 once we
know that the action of an iwip $\phi\in\On$ on $\partial_\infty(X_n)$
has north-south dynamics. And this was proved by Levitt and Lustig
\cite{LeL--03}. Thus we have:

\proclaim{2.6\ Theorem}
If $n \ge 3$, then
$\text{\rm{Out}}(F_n)$ is a Powers group;
in particular its reduced $C^*$-algebra is simple.
\endproclaim
 
If $\Gamma\subset\On$ is torsion-free, the stabilizers in $\Gamma$ of the
fixed points of iwips
are cyclic, so the arguments above establish:

\proclaim{2.7\ Theorem}
If $\Gamma\subseteq\On$ is  torsion-free, non-cyclic and contains an iwip,
then it is  a Powers group;
in particular its reduced $C^*$-algebra is simple.
\endproclaim

 Let $\Gamma$ be a group, $A$ a unital $C^*$-algebra,
and $\alpha$ an action
of $\Gamma$ on $A$ such that $A$ contains no non-trivial
$\alpha(\Gamma)$-invariant two-sided
ideals.
In general, the reduced crossed product
$A \rtimes_{\alpha,r} \Gamma$ need not be simple.
However, it  will always be simple
if $\Gamma$ is a Powers group \cite{HaS--86}.
Thus we have the following corollary.

\proclaim{2.8\ Corollary} Let $\Gamma$ be a group as in Theorem 2.2, 2.3,
2.6 or 2.7. If $\Gamma$ acts on a unital
$C^*$-algebra $A$ and leaves no non-trivial 2-sided ideals
invariant, then the corresponding reduced crossed
product is a simple $C^*$-algebra.
\endproclaim

This corollary applies, for example, to the reduced crossed
product $\Cal C (L_{\Gamma}) \rtimes_r \Gamma$ associated to the
action of $\Gamma$ on the algebra of continuous functions
on the limit set of $\Gamma$ in $\Cal{PMF}\!_S$ or $\partial_\infty(X_n)$.
\medskip

   {\it Remarks on minimality and compactness.}
 In previous papers, Condition (ii) of Proposition 1.1
is replaced by other conditions involving minimality for the action of
$\Gamma$ on $\Omega$ and a compactness assumption on $\Omega$.
As we have seen, compactness is not necessary.
Minimality can always be obtained, if desired, by replacing
$\Omega$ with the closure of the set of fixed points of hyperbolic
elements. We remark that the action of $\Mod_S$ on
$\PMF\!_S$ is minimal for all $g \ge 1$ (see e.g. \S \ VII of
expos\'e 6
in \cite{FLP--79}), but the action of $\On$ on $\partial_\infty(X_n)$
is not: the ideal points of the simplicial spine
$K_n\subset X_n$ defined in \cite{CuV--86} form a closed invariant
subset, for example.

\medskip

   {\it The Limitations of the Method.} It would be
interesting
to understand more precisely the
limitations of the elementary method used
here to establish the
$C^*$-simplicity of groups. We used Corollary 1.2 in a rather
weak form: in our examples the stabilizers of the endpoints
of our hyperbolic
homeomorphisms were virtually cyclic, whereas Corollary 1.2
would allow large malnormal
\footnote{
A subgroup $H$ of a group $G$ is {\it malnormal}
if $gHg^{-1} \cap H = \{e\}$ for any $g \in G$ such that $g \notin H$.
}
stabilizers, for example. Since
the centralizer of a hyperbolic element must stabilize its
fixed points, a natural challenge arises in the case of $F_n\times F_n$,
where all centralizers contain a copy of $\Bbb Z^2$ but are not
malnormal. Does this group admit a faithful action satisfying the
conditions of Proposition 1.1?
(Note that the direct product of two non-abelian free groups
is $C^*$-simple since a spatial tensor product of simple $C^*$-algebras
is a simple $C^*$-algebra \cite{Tak--64}.)

Other challenging examples are the groups
$\text{\rm{PSL}}(n,\Bbb Z)$ for $n \ge 3$:
we know that these groups are $C^*$-simple (by \cite{BCH--94}),
but we do not know whether they are Powers groups.

\bigskip
\head
3. Appendix on the sufficiency of Powers' criterion
\endhead
\medskip

For the convenience of the reader, we include a  proof of the result of
Powers  quoted at the
beginning of
Section~1. Powers' original proof was formulated only for non-abelian free
groups
\cite{Pow--75}, but the following adaptation is entirely straightforward.

\proclaim{3.1\ Theorem} Powers groups are $C^*$-simple.
\endproclaim

\demo{Proof (following  \cite{Pow--75})} Consider a Powers group
$\Gamma$, its left-regular representation $\lambda_{\Gamma}$,
 a non-zero ideal
$\Cal I$ of its reduced $C^*$-algebra,
 and an element $U \ne 0$ in $\Cal
I$. We want
to show  that $\Cal I$ contains an element
 $Z$ such that $\norm{Z - 1} <
1$, and in
particular such that $Z$ is invertible
 (with inverse $\sum_{n=0}^{\infty}Z^n$).
\par

   Upon replacing $U$ by a scalar multiple of $U^*U$,
we may assume that
$U = 1 + X$ and $X = \sum_{x \in \Gamma, x\neq e}
z_x \lambda_{\Gamma}(x)$,
with $z_x \in \Bbb C$.
Choose $\epsilon,\delta$ with $0 < \epsilon < \delta \le 1$
(here, we could set $\delta = 1$, but we will use
the freeness in choosing $\delta$ in the next proof).
Then there exists a finite subset  $F$ of
$\Gamma \smallsetminus \{e\}$ such that, if
$$
   X' \, = \, \sum_{f \in F} z_f \lambda_{\Gamma} (f) ,
$$
then $\norm{X' - X} < \epsilon$.
Set $U' = 1 + X'$, so that $\norm{ U' - U} < \epsilon$.
Choose an integer $N$ so large that
$\frac{2}{\sqrt N} \norm{X'} < \delta  - \epsilon$.
\par

   Let now $\Gamma = C \sqcup D$ and $\gamma_1,\hdots,\gamma_N$
be as in the definition of a Powers group. Set
$$
\aligned
 V \, = \,  \frac{1}{N}
    \sum_{j=1}^N \lambda_{\Gamma}(\gamma_j) U
\lambda_{\Gamma}(\gamma_j^{-1})
    \qquad\qquad &
 V' \, = \, \frac{1}{N}
    \sum_{j=1}^N \lambda_{\Gamma}(\gamma_j) U'
\lambda_{\Gamma}(\gamma_j^{-1})
 \\
 Y \, = \,  \frac{1}{N}
    \sum_{j=1}^N \lambda_{\Gamma}(\gamma_j) X
\lambda_{\Gamma}(\gamma_j^{-1})
 \qquad\qquad &
 Y' \, = \, \frac{1}{N}
    \sum_{j=1}^N \lambda_{\Gamma}(\gamma_j) X'
    \lambda_{\Gamma}(\gamma_j^{-1}).
\endaligned
$$
Note that $V = 1+Y \in \Cal I$ and $V' = 1 + Y'$. We show below that
$\norm{Y'} < \delta -  \epsilon$. This implies that
$\norm{Y} \le \norm{Y'} + \norm{Y-Y'} \le \norm{Y'} + \norm{X-X'}
< \delta \le 1$.
As $\Cal I$ contains the invertible element $V = 1 + Y$,
the $C^*$-algebra $C^*_{\lambda}(\Gamma)$ is indeed simple.  \par

 For $j \in \{1,\hdots,N\}$, denote by $P_j$ the orthogonal projection of
$\ell ^2 (\Gamma)$ onto $\ell ^2 (\gamma_j D)$.
We have
$$
   (1-P_j) \lambda_{\Gamma}(\gamma_j) X' \lambda_{\Gamma}(\gamma_j^{-1})
(1-P_j)
   \, = \, 0 ;
$$
indeed, since $fC \cap C = \emptyset$ for all $f \in F$, we have
$$
\aligned
 \bigg( \lambda_{\Gamma}(\gamma_j)X'
    \lambda_{\Gamma}(\gamma_j^{-1})(1-P_j) \bigg) (\ell^2(\Gamma))
   \ &\subset \
   \bigg( \lambda_{\Gamma}(\gamma_j)X' \bigg) (\ell^2(C)) \\
   \ &\subset \
   \bigg( \lambda_{\Gamma}(\gamma_j) \bigg) (\ell^2(D) ) \ = \ P_j
(\ell^2(\Gamma)) .
\endaligned
$$
It follows that
$$
  V' \, = \, 1 \, + \,
    \frac{1}{N} \sum_{j=1}^N P_j
    \lambda_{\Gamma}(\gamma_j) X' \lambda_{\Gamma}(\gamma_j^{-1})
    \, + \, \left(
    \frac{1}{N} \sum_{j=1}^N P_j
    \lambda_{\Gamma}(\gamma_j) X' \lambda_{\Gamma}(\gamma_j^{-1})
              (1-P_j) \right) ^* .
$$
Since the subsets $\gamma_jD$ of $\Gamma$ are pairwise disjoint, the
operators
$X'_j \Doteq P_j \lambda_{\Gamma}(\gamma_j) X'
\lambda_{\Gamma}(\gamma_j^{-1})$
have pairwise orthogonal ranges in $\ell ^2 (\Gamma)$, and we have
$$
\norm{ \frac{1}{N} \sum_{j=1}^N X'_j }
   \, \le \,
\frac{1}{ \sqrt N } \max_{j=1}^n \norm{ X'_j }
   \, \le \,
\frac{1}{ \sqrt N } \norm{ X' } .
$$
Similarly
$$
\norm{ \left( \frac{1}{N} \sum_{j=1}^N X'_j (1-P_j) \right) ^* }
   \, = \,
\norm{ \frac{1}{N} \sum_{j=1}^N X'_j (1-P_j) }
   \, \le \,
\frac{1}{\sqrt N}\norm{X'}.
$$
Consequently
$$
\norm{Y'}
   \, = \,
\norm{ V' - 1 }
   \, \le \,
\frac{2}{ \sqrt N } \norm{X'}
   \, < \,
\delta - \epsilon .
$$
As already observed,  this completes the proof. $\square$
\enddemo

In Section 6 of \cite{BekL--00} the above
proof  is recast in the language of functions of positive type.

\bigskip

   A linear form $\tau$ on $C^*_\lambda(\Gamma)$ is a {\it normalised trace}
if $\tau(1)=1$ and $\tau(U^*U)\ge 0,\ \tau(UV)=\tau(VU)$
for all $U,V \in C^*_\lambda(\Gamma)$.
We have $\abs{\tau(X)} \le \norm{X}$ for all $X \in C^*_\lambda(\Gamma)$
(see, e.g., Proposition 2.1.4 in \cite{D$C^*$--69}).
The {\it canonical trace} is uniquely defined by
$$
\tau_{\text{\rm{can}}}\left(\sum_{f\in F}z_f\lambda_\Gamma(f)\right)
\, = \, z_e
$$
for every {\it finite} sum $\sum_{f\in F}z_f\lambda_\Gamma(f) $ where
$z_f\in\Bbb C$ and $F \subset \Gamma$ contains $e$.

\proclaim{3.2\ Proposition} If $\Gamma$ is a Powers group, then the 
canonical
trace
is the only  normalized trace on $C^*_\lambda(\Gamma)$.
\endproclaim

\demo{Proof} Let $\tau$ be a normalised trace on $C^*_\lambda(\Gamma)$
and let $X \in C^*_\lambda(\Gamma)$ be such that
$\tau_{\text{\rm{can}}}(X) = 0$. It is enough to show that $\tau(X) = 0$.

   Choose $\delta > 0$. The previous proof shows that there exist
$N \ge 1$ and $\gamma_1,\hdots,\gamma_N \in \Gamma$
such that $\norm{Y} \le \delta$ for
$$
Y \, = \,
\frac{1}{N} \sum_{j=1}^N \lambda_{\Gamma}(\gamma_j) X
\lambda_{\Gamma}(\gamma_j^{-1}) .
$$
As $\tau(Y) = \tau(X)$, we have $\abs{\tau(X)} \le \norm{Y} < \delta$.
As $\delta$ is arbitrary, this implies $\tau(X) = 0$.
$\square$
\enddemo

\enddocument